\documentclass[12pt]{amsart}

\usepackage{latexsym}
\usepackage{amsfonts}
\usepackage{amsmath}
\newtheorem{theorem}{Theorem}[section]
\newtheorem{rem} [theorem] {Remark}
\newtheorem{prop} [theorem]{Proposition}
\newtheorem{lemma}[theorem]{Lemma}

\newtheorem{definition}[theorem]{Definition}
\newcommand{\ovprt}{\overline{\partial}}
\newcommand{\ovli}{\overline}

\numberwithin{equation}{section}
\begin{document}
\title{Compactness of the solution operator to $\ovprt $ in weighted $L^2$ - spaces.}

\author{Friedrich Haslinger and Bernard Helffer}

\thanks{Both authors are grateful to the Erwin Schr\"odinger Institute, Vienna, where most of this work was carried out.}

\address{F. Haslinger: Institut f\"ur Mathematik, Universit\"at Wien,
Nordbergstrasse 15, A-1090 Wien, Austria}
\email{friedrich.haslinger@univie.ac.at}

\address{B. Helffer : Laboratoire de Math\'ematique, Universit\'e Paris Sud, B\^atiment 425, F-91405 Orsay, France}
\email{bernard.helffer@math.u-psud.fr}


\maketitle

\begin{abstract} ~\\
In this paper we discuss compactness of the canonical solution operator
to $\ovprt $ on weigthed $L^2$ spaces on $\mathbb C^n.$
For this purpose we apply ideas which were used for the Witten Laplacian in the real case and various methods of spectral theory of these operators. We also point out connections to the theory of Dirac and Pauli operators.

\end{abstract}

\vskip 1 cm

\section{Introduction.}~\\

\subsection{Background for bounded pseudoconvex domains}~\\

Let $\Omega $ be a bounded pseudoconvex domain in $\mathbb C^n.$ 
We consider the 
$\ovprt $-complex 
$$
L^2(\Omega )\overset{\ovprt }\longrightarrow L^2_{(0,1)}(\Omega)
\overset{\ovprt }\longrightarrow \dots \overset{\ovprt }\longrightarrow
L^2_{(0,n)}(\Omega)\overset{\ovprt }\longrightarrow 0\, ,  
$$
where $L^2_{(0,q)}(\Omega)$ denotes the space of $(0,q)$-forms on $\Omega$ with
coefficients in $L^2(\Omega).$ The $\ovprt $-operator on $(0,q)$-forms is given by
$$
\ovprt \left ( \sum_J\,^{'} a_J  \, d\ovli z_J \right )\sum_{j=1}^n \sum_J\,^{'}\  \frac{\partial a_J}{\partial \ovli z_j}d\ovli z_j\wedge
d\ovli z_J,$$
where $\sum \ ^{'} $ means that the sum is taken only over increasing multi-indices $J.$

The derivatives are taken in the sense of distributions, and the domain of $\ovprt $
consists of those $(0,q)$-forms for which the right hand side belongs to
$L^2_{(0,q+1)}(\Omega).$ Then $\ovprt $ is a densely defined closed operator, and
therefore has an adjoint operator from $L^2_{(0,q+1)}(\Omega)$ into
$L^2_{(0,q)}(\Omega)$ denoted by $\ovprt ^* .$

The complex Laplacian 
$\Box = \ovprt \, \ovprt ^* + \ovprt ^* \, \ovprt $
acts as
an unbounded selfadjoint  operator on 
$L^2_{(0,q)}(\Omega ),\ 1\le q \le n,$
it is surjective and therefore
has a continuous inverse, the $\ovprt $-Neumann operator $N_q.$ If $v$ is a closed
$(0,q+1)$-form, then $\ovprt ^* \, N_{q+1}v$ provides the canonical solution to 
$\ovprt u = v,$ which is orthogonal to the kernel of $\ovprt$ and
so has minimal norm (see for instance \cite{ChSh}).

\vskip 0.5 cm

A survey of the $L^{2}$-Sobolev theory
of the $\overline{\partial}$-Neumann problem is given in \cite{BS}.
\vskip 0.5 cm

The question of compactness of $N_q$ is of interest for various
reasons. For example, compactness of $N_q$ implies global regularity
in the sense of preservation of Sobolev spaces \cite{KN}. Also, the
Fredholm theory of Toeplitz operators is an immediate consequence
of compactness in the $\overline{\partial}$-Neumann problem \cite{V},
\cite{HI}, \cite{CD}. There are additional ramifications for certain
$C^{*}$-algebras naturally associated to a domain in $\mathbb{C}^{n}$
\cite{SSU}. Finally, compactness is a more robust property than
global regularity - for example, it localizes, whereas global
regularity does not - and it is generally believed to be more
tractable than global regularity. 

A thorough discussion of
compactness in the $\overline{\partial}$-Neu\-mann problem can be found
in \cite{FS2}.

The study of the $\overline{\partial}$-Neumann problem is essentially
equivalent  to the study of the
canonical solution operator to $\overline{\partial}$:

The $\ovprt $-Neumann operator $N_q$ is compact from $L^2_{(0,q)}(\Omega )$ to itself if and only if the canonical solution operators
$$\ovprt ^*N_q: L^2_{(0,q)}(\Omega )\longrightarrow L^2_{(0,q-1)}(\Omega )   \ \ \textrm{and} \ \
\ovprt ^*N_{q+1}: L^2_{(0,q+1)}(\Omega )\longrightarrow L^2_{(0,q)}(\Omega )$$
are compact.

\medskip

 Interestingly,
in many situations, the restriction of the canonical solution operator
to forms with \emph{holomorphic} coefficients arises naturally
\cite{SSU}, \cite{FS1}. 
Compactness of the restriction to forms with holomorphic coefficients already implies compactness
of the original solution operator to $\ovprt $ in the case of convex domains,
see \cite{FS2}. There are many examples for non-compactness, where the obstruction 
already occurs for forms with holomorphic coefficients (see \cite{Has1}, \cite{Has2},  \cite{Kr} and \cite{L}). 

\medskip

In \cite{CD} it is shown that compactness of the $\ovprt $-Neumann operator implies compactness of the commutator $[P,M],$ where $P$ is the Bergman projection and $M$ is pseudodifferential operator of order $0.$
In \cite{Has4} it is shown that compactness of the canonical solution operator to $\overline \partial $ restricted to $(0,1)$-forms with holomorphic coefficients implies compactness of the commutator $[P,M]$ defined on the whole $L^2(\Omega ).$

Let $A_{(0,1)}^2(\Omega ) $ denote the space of 
all $(0,1)$-forms with holomorphic coefficients
belonging to $L^2(\Omega ).$

Another result from \cite{Has4} states that  for a bounded pseudoconvex domain $\Omega $ the $\ovprt $-Neumann operator
$$N_1: L^2_{(0,1)}(\Omega )\longrightarrow L^2_{(0,1)}(\Omega ) $$
restricted to $(0,1)$-forms with holomorphic coefficients can be written in the form
$$\mathcal{P}N\mathcal{P}f = \sum_{k=1}^n [P,M_k] \left( \sum_{j=1}^n [\overline M_j,P]f_j\right) d\overline z_k 
$$
here $\mathcal{P}:  L^2_{(0,1)}(\Omega )\longrightarrow A^2_{(0,1)}(\Omega )$ denotes the componentwise projection and 
$M_j$ and $\overline M_j$ denotes the multiplication by $z_j$ and $\overline z_j$ respectively. 

\medskip

The restriction of the canonical solution operator to forms with holomorphic 
coefficients has many interesting aspects, which in most cases correspond to 
certain growth properties of the Bergman kernel. 
\medskip 

In \cite{Has1} the canonical solution operator $S_1$ to $\ovprt $
restricted to $(0,1)$-forms with holomorphic coefficients is investigated. 

It is shown that
the canonical solution operator
$S_1 : A_{(0,1)}^2(\Omega ) \longrightarrow L^2(\Omega )$
has the form
$$
S_1(g) (z) = \int_{\Omega }B(z,w)<g(w), z-w>\,d\lambda (w),$$
where $ B$ denotes the Bergman kernel of $\Omega $ and
$$<g(w), z-w>=\sum_{j=1}^n g_j(w) (\overline z_j -\overline w_j),$$
for $z=(z_1,\dots ,z_n)$ and $w=(w_1,\dots ,w_n);$
it can also be written in the form
$$S_1g= \sum_{j=1}^n [\overline M_j,P]g_j .$$
It follows that the canonical
solution operator is a Hilbert Schmidt operator
for the unit disc $\mathbb D $ in $\mathbb C ,$ 
but fails to be Hilbert Schmidt for the unit 
ball in $\mathbb C^n, \ n\ge 2$ (see also \cite{LY}).

\vskip 0.5 cm
\medskip
\subsection{The case of unbounded domains}~\\
\vskip 0.3 cm

Not very much is known in the case of unbounded domains.
 In this paper we discuss the compactness of the canonical solution operator to $\ovprt $ on weighted
$L^2$-spaces over $\mathbb C^n.$ We define
$$L^2(\mathbb C^n, \varphi) = \{ f:\mathbb C^n \longrightarrow \mathbb C \ : \int_{\mathbb C^n} \vert f(z) \vert^2 \exp(-2\varphi (z))\, d\lambda (z) < \infty \}, $$
where $\varphi $ is a suitable weight-function.

There is an interesting  connection of $\ovprt $ with the theory of Schr\"odinger operators with magnetic fields,
see for example \cite{Ch}, \cite{B}, \cite{FS3} and \cite{ChF} for recent contributions exploiting this point of view.

For the case of one complex variable  results of Helffer-Mohamed \cite{HeMo}, Iwatsuka
\cite{I}  and Shen \cite{She1} can be used  to discuss compactness of the canonical solution operator to $\ovprt.$
For instance, if $\varphi (z)=|z|^2,$ then the canonical solution operator
$S:L^2(\mathbb C, \varphi )\longrightarrow L^2(\mathbb C, \varphi )$ to $\ovprt $
fails to be compact.
If  $\Delta \varphi (z)\to \infty$
as $|z|\to \infty ,$ then the
canonical solution operator $S\,:\,L^2(\mathbb C, \varphi )\longrightarrow L^2(\mathbb C, \varphi )$ to $\ovprt $ is
compact (\cite{Has3} ). 

In this paper
we first give a necessary and sufficient condition in terms of the weight function $\varphi $ in the complex one-dimensional case for the  solution operator to be compact on $L^2(\mathbb C, \varphi)$ continuing the work from \cite{Has3} and using results from
\cite{AuBe}, \cite{HeMo}, \cite{I}, \cite{She1}, and \cite{St}.

In the case of several complex variables, we meet an obvious condition
for solving $\ovli\partial u=f$. The $(0,1)$-form $f$ should satisfy
$\ovli \partial f =0$. So we are asking for the existence of
a continuous operator $S^{can}$, which will be called the canonical
solution operator~:
\begin{equation}\label{defcan}
L^2_{(0,1)}(\mathbb C^n,\varphi) \cap \mbox{ Ker }  \ovli\partial
 \ni f \mapsto u=S^{can} f \in L^2(\mathbb C^n, \varphi) \cap ( \mbox{
 Ker } 
 \ovli \partial)^{\perp}\;,
\end{equation}
giving the minimal solution of the problem.\\
When the weight function $\varphi $ is  plurisubharmonic, we will for
example  show that the condition that the lowest eigenvalue $\lambda_\varphi $ of the Levi matrix $M_\varphi $ satisfies
$$\lim_{|z|\to \infty } \lambda_\varphi (z) = +\infty$$
implies the existence of  the canonical solution operator and its compactness.

For  decoupled weights
$$
\varphi (z)=\varphi_1(z_1) + \varphi_2(z_2) + \dots + \varphi_n(z_n)
\,,
$$ 
the canonical solution operator to $\ovprt $ fails, under very weak
additional assumptions to be compact and we will show that it is  even
true on $A^2_{(0,1)}(\mathbb C^n, \varphi )$ (see \cite{Sch}).\\

There are other interesting connections between $\overline{\partial}$ and
Schr\"{o}d\-inger operators, see for example the
discussion in \cite{B} and  between
compactness in the $\overline{\partial}$-Neumann problem and
property (P) on the one hand, and the asymptotic behavior, in a
semi-classical limit, of the lowest eigenvalues of certain magnetic
Schr\"{o}d\-inger operators and of their non-magnetic counterparts,
respectively, on the other (\cite{FS3}). The
main result in \cite{FS3} shows that (for certain Hartogs domains in
$\mathbb{C}^{2}$) 
compactness properties of the $\ovprt $-Neumann operator may be interpreted as a consequence of well known
diamagnetic inequalities (originally due to Kato) 
 in the theory of Schr\"{o}dinger operators (see \cite{LL},
 \cite{CFKS} and 
\cite{He1}).
\vskip 0.2 cm

Finally, we also point out some interesting connections to the theory
of Dirac and Pauli operators, when discussing the case of non-compact
resolvents (see \cite{CFKS}, \cite{Er}, \cite{HNW}, \cite{Roz},
\cite{Tha}).\\

{\bf Acknowledgements}~\\
The authors would like to thank P.~Auscher and Z.~Shen for useful discussions.

\vskip 1 cm
\section{The complex one-dimensional case}

\vskip 1 cm 

Let $\varphi $ be a subharmonic $\mathcal C^2$-function.
We want to solve $\ovprt u=f$ for $f\in L^2(\mathbb C, \varphi ).$ The canonical solution
operator to $ \ovprt $ gives a solution with minimal $L^2(\mathbb C, \varphi )$-norm.
We substitute $v=u\,e^{-\varphi }$ and $g=f\,e^{-\varphi }$
and the equation becomes
$$\overline D v=g \ ,$$
where
\begin{equation}
 \overline D = e^{-\varphi }\, \frac {\partial }
{\partial \overline z}\, e^{\varphi }.
\end{equation}
$u$ is the minimal solution to the $\ovprt $-equation in $L^2(\mathbb C, \varphi )$ if and only if
$v$ is the solution to $\overline D v=g$ which is minimal in $L^2(\mathbb C )\,$.

The formal adjoint of $\overline D$ is 
\begin{equation}
D=-e^{\varphi }\frac{\partial}{\partial z}e^{
-\varphi }.
\end{equation} 
Let us introduce
\begin{equation}\label{defS}
\mathcal S =  \overline D D\;.
\end{equation}

 Since $\overline D=\frac{\partial }{\partial \overline z}+
\frac{\partial \varphi }{\partial \overline z}$ and $D=-\frac{\partial }{\partial z}+
\frac{\partial \varphi }{\partial z}$,  we see that
$$ \mathcal S =-\frac{\partial^2}{\partial z \partial \overline z}-\frac{\partial
\varphi }{\partial \overline z}\,\frac{\partial}{\partial z}+ \frac{\partial \varphi }
{\partial z}\,\frac{\partial}{\partial \overline z}+\left |\frac{\partial \varphi }
{\partial z} \right |^2+\frac{\partial^2 \varphi}{\partial z \partial
\overline z}\;.$$
So
\begin{equation}\label{formuleS}
\mathcal S =\frac{1}{4}\, (-\Delta_A + B ),
\end{equation}
where the $1$-form $A=A_1\,dx+A_2\,dy$ is related to the weight
$\varphi$ by
\begin{equation}
A_1=-\partial_y \varphi \,,\, A_2 = \partial_x   \varphi \,,
\end{equation} 
\begin{equation}
\Delta_A =  \left ( \frac{\partial}{\partial x}-i A_2 \right )^2 +
\left ( \frac{\partial}{\partial y}+i A_1 \right )^2 \;,
\end{equation}
and the magnetic field $B dx\wedge dy$ satisfies
\begin{equation}
B(x,y) =\Delta \varphi (x,y)\;.
\end{equation}

 Hence $\mathcal S $ is (up to a multiplicative constant) a Schr\"odinger
operator with magnetic field and  an electric potential $B$. In
addition, we know from \cite{Sima} that this operator is essentially
 self-adjoint on $C_0^\infty(\mathbb C)$.\\

In \cite{Has3} (completing a result of M.~Christ \cite{Ch}), a link
was established between the compactness  of the canonical solution
operator to $\ovprt $ and the properties of the resolvent of $\mathcal S$. In this setting it was supposed that the weight
functions $\varphi $
 are in the class $\mathcal W$.
\begin{definition}\label{defW}~\\
We say that $\varphi$ is in the class $\mathcal W$ if~:
\begin{enumerate}
\item   $\nu =\Delta  \varphi
dx ,$ is a doubling measure, which means that
there exists a constant $C$ such that for all $z\in \mathbb C$ and $r\in \mathbb R^+,$
$$\nu (B(z,2r))\le C\; \nu (B(z,r)),$$
where $B(z,r)$ denotes the ball with center $z$ and radius $r$\;;
\item there exists
a constant $\delta >0$ such that for all $z\in \mathbb C,$
$$\nu (B(z,1))\ge \delta.$$
\end{enumerate}
\end{definition}
In fact Marco, Massaneda and Ortega-Cerda \cite{MMO} (Theorem C, p.~884) found out that already condition (1) in the last definition implies that the canonical solution operator to $\ovprt $
is continuous. Hence it follows from \cite{Has3}

\begin{theorem}\label{thhaschr} ~\\
Let $\varphi $ be a subharmonic $\mathcal C^2$-function on $\mathbb R^2$ such that $\Delta \varphi $ defines a doubling measure. The canonical solution operator
$S:L^2(\mathbb C ,\varphi )\longrightarrow L^2(\mathbb C, \varphi )$
to $\ovprt $ is compact if and only if $\mathcal S$ has compact resolvent.
\end{theorem}

Now we prove a criterion of compactness, which can be expressed  in terms of the weight function $\varphi $ only.
Here we extend a result due to Helffer and Morame (\cite{HeMor}) based on methods developed by Iwatsuka (\cite{I}) and Shen (\cite{She1}).
For this purpose we assume the stronger condition that the weight
function $\varphi $ is a subharmonic $\mathcal C^2$ function and that $\Delta\varphi $ belongs
to the reverse H\"older class $B_2(\mathbb R^2) $ consisting of $L^2 $
positive and almost non zero everywhere functions $V$ for which there exists a constant $C>0$ such that
$$\left( \frac{1}{|Q|} \int _Q V^2\,dx\right)^{\frac 12} \leq
 C \left(\frac {1}{|Q|} \int_Q V \, dx\right)$$
for any ball $Q$ in $\mathbb R^2$\,.

It is known that if $V$ is in $B_q$ for some $q>1$
 then $V$ is in the Muckenhoupt class $A_\infty$ and the corresponding measure
 $V(x)dx$
 is doubling. More precisely it is known from \cite{St} that
$$
A_\infty =\cup_{q>1} B_q\;.
$$
 Note that any positive (non zero) polynomial is in $B_q$ for any
 $q>1$.

\begin{theorem}\label{hasimp}~\\
Let $\varphi $ be a subharmonic $\mathcal C^2$- function on 
$\mathbb R^2$ such that  
\begin{equation}\label{assB2}
\Delta\varphi \in B_2(\mathbb R^2)\;.
\end{equation}
Then the canonical solution operator
$S:L^2(\mathbb C ,\varphi )\longrightarrow L^2(\mathbb C, \varphi )$ to $\ovprt $ is compact if and only if
  \begin{equation}\label{2.9}
\lim_{|z|\rightarrow \infty} \int_{B(z,1)}
\Delta\varphi (y)\, dy = +\infty \,.
\end{equation}
\end{theorem}

\begin{proof}~\\
Using Theorem \ref{thhaschr}, we have just to analyze if  $-\Delta_A +
 \Delta\varphi$
 has compact resolvent.\\
Using the standard comparison between selfadjoint operators~:
\begin{equation} \label{comparison}
- 2 \Delta_A \geq -\Delta_A + \Delta\varphi \ge -\Delta_A \,
\end{equation}
we observe that $-\Delta_A + \Delta \varphi$ has  compact resolvent
 if and only if $-\Delta_A$ has compact resolvent.\\

In one direction,  we can apply a result of Iwatsuka (\cite{I} ,
Theorem 5.2) which says 
\begin{prop}\label{neccond}~\\
Suppose that $A\in H^1_{loc}$ and that  $-\Delta_A$ has compact resolvent.
Then
\begin{equation}\label{nec}
\lim_{|z|\rightarrow  \infty} \int_{B(z,1)}B(y)^2 \,dy = +\infty\;,
\end{equation}
with  $B=\text{curl} A$.
\end{prop}
 Iwatsuka adds a $C^\infty$ assumption on the magnetic potential. But
 at least in the two dimensional case, one can use properties of the Curl operator as mentioned in Appendix I
of \cite{T}, in order to release this assumption. Note that  in our
 case $B=\Delta \varphi$.
 By the definition of
 the reverse H\"older class $B_2(\mathbb R^2)$, \eqref{nec}  implies \eqref{2.9}.

\vskip 0.2 cm
For the other direction, we first use a version of the diamagnetic property for
Schr\"odinger operators (see for example \cite{KS} Cor. 1.4) saying that~:\\
If $-\Delta + \Delta \varphi$ has compact resolvent, then
 $-\Delta_A + \Delta \varphi$ has compact resolvent.
 So it is enough to prove that $-\Delta + V$
 has compact resolvent with $V=\Delta \varphi$.\\

By the Main Theorem in \cite{I},  it suffices to show that
\begin{equation}\label{condb}
\lim_{|z|\rightarrow \infty}  \lambda_{0,V} (B(z,1))= +\infty\;,
\end{equation}
 where  $\lambda_{0,V}(B(z,1))$ is
 the lowest eigenvalue of the Dirichlet realization of $-\Delta + V$
 in $B(z,1)\,$.
Without loss of generality, we can consider, instead of  balls, cubes. In
 this case  we use the following improved 
version of the Fefferman-Phong Lemma as given in  \cite{AuBe}.
 \begin{lemma}\label{LemmaAuscher}~\\
 If $V \in A_\infty$, then there  exists $C_V>0$ and $\beta_V \in ]0,1[$
  such that, for all cubes  $Q$ (with sidelength $R$), for all $u\in C_0^\infty(Q)$,
\begin{equation}\label{ineqfond}
 C_V\frac{m_\beta (R^2 \Theta_Q)}{R^2}   \int|u(y)|^2\, dy \leq
  \int (|\nabla u(y)|^2 + V (y) |u(y)|^2)\;dy
 \end{equation}
where 
$$
\Theta_Q= \frac{1}{|Q|} \int_Q  V(y) dy \;,
$$
and 
$$
m_\beta(t) = t\mbox{ for } t\leq 1\,,\; \mbox{ and } m_\beta(t) = t^{\beta_{V}}
 \mbox{ for } t \geq 1\;.
$$
\end{lemma}
We apply Lemma \ref{LemmaAuscher} with $R=1$ and $V=\Delta
 \varphi$. \eqref{ineqfond} gives a lower bound for $\lambda_{0,V}
 (Q)$ by $C_V \Theta_Q^{\beta_V}$ each time that $\Theta_Q \geq 1$. 
 Therefore  Assumption \eqref{2.9}  implies \eqref{condb}  and we are done.
\end{proof}

\begin{rem}~\\
As a variant of the proof, we have the following statement. Suppose that $\Delta
\varphi$
 belongs to $A_\infty$ (at $\infty$, i.e. for all the balls meeting
 the complement of a compact $K$) and that
 $$\liminf
 _{|z|\rightarrow \infty} \int_{B(z,1)} \Delta \varphi (y) dy >0\;,$$
then the canonical solution operator $S$ is well defined
 and Theorem \ref{thhaschr} is true.\\
Note that we have also shown that if  $$\lim
 _{|z|\rightarrow \infty} \int_{B(z,1)} \Delta \varphi (y) \, dy =  +\infty\;,$$
then $S$ is compact.\\
 
We learn from Z.~Shen, that, in this $2$-dimensional case, one can, by
other techniques developed in \cite{She2}, improve the necessary part
 due to Iwatsuka and deduce the same result under the weaker assumption that $\Delta
 \varphi \in A^\infty$. This proof is much more involved
 and strongly limited to the two-dimensional case.

\end{rem}

\vskip 1 cm

\section{The $\ovprt $-equation  in weighted $L^2$ - spaces of several
  complex variables~: the  canonical solution operator.}

\vskip 1 cm
  
Here we apply ideas which were used in the analysis of Witten Laplacian in the real case, see \cite{HeNi}.

Let $\varphi : \mathbb C^n \longrightarrow \mathbb R $ be a $\mathcal C^2$-weight function and define the space
$$L^2(\mathbb C^n , \varphi )=\{ f:\mathbb C^n \longrightarrow \mathbb C \ : \ \int_{\mathbb C^n}
|f|^2\, e^{-2\varphi}\,d\lambda < \infty \},$$
the space $L^2_{(0,1)}(\mathbb C^n, \varphi )$ of $(0,1)$-forms with coefficients in
$L^2(\mathbb C^n , \varphi )$ and the space $L^2_{(0,2)}(\mathbb C^n, \varphi )$ of $(0,2)$-forms with coefficients in
$L^2(\mathbb C^n , \varphi ).$

Let $A^2(\mathbb C^n , \varphi )$ denote the space of entire functions belonging to $L^2(\mathbb C^n , \varphi )$.

We consider the 
$\ovprt $-complex 
$$
L^2(\mathbb C^n , \varphi )\overset{\ovprt }\longrightarrow L^2_{(0,1)}(\mathbb C^n , \varphi )
\overset{\ovprt }\longrightarrow L^2_{(0,2)}(\mathbb C^n,\varphi) \,.
$$

For $v\in L^2(\mathbb C^n )$, let
$$\overline D_1 v = \sum_{k=1}^n \left( \frac{\partial v}{\partial \ovli z_k}+
\frac{\partial \varphi}{\partial \ovli z_k}\, v \right) \, d\ovli z_k$$
and for $g=\sum_{j=1}^n g_j\, d\ovli z_j \in L^2_{(0,1)}(\mathbb C^n ) $, let
$$\ovli D^*_1 g = \sum_{j=1}^n \left( \frac{\partial \varphi}{\partial z_j}\, g_j
-\frac{\partial g_j}{\partial z_j} \right) ,$$
where the derivatives are taken in the sense of distributions.\\
It is easy to see that $\ovprt u =f$ for $u\in L^2(\mathbb C^n , \varphi )$ and
$f\in L^2_{(0,1)}(\mathbb C^n, \varphi )$ if and only if $\ovli D_1 v = g, $ where
$v= u\, e^{-\varphi }$ and $g= f\, e^{-\varphi }.$ It is also clear that
the necessary condition $\ovprt f=0$ for solvability holds if and only if
$\ovli D_2 g =0 $ holds. Here
$$\ovli D_2 g = \sum_{j,k=1}^n \left ( \frac{\partial g_j}{\partial \ovli z_k}
+\frac{\partial \varphi }{\partial \ovli z_k}\,g_j \right ) \, d\ovli z_k \wedge
d\ovli z_j.$$
So the existence and the  analysis of the canonical solution operator introduced in
\eqref{defcan} is equivalent to the existence and the analysis
 of the canonical solution operator for  $\ovli D$, the equivalence being given by
\begin{equation}\label{equivscan}
S^{can}_\varphi = \exp (- \varphi ) \,  S^{can} \, \exp (\varphi )\;.
\end{equation}

We consider the corresponding $\ovli D$-complex with in particular~:
$$ 
L^2(\mathbb C^n ) \underset{\underset{\ovli D_1^* }{\longleftarrow }}{\overset{\ovli D_1 }{\longrightarrow}}
L^2_{(0,1)}(\mathbb C^n )   \underset{\underset{\ovli D_2^* }{\longleftarrow }}{\overset{\ovli D_2 }{\longrightarrow}}
L^2_{(0,2)}(\mathbb C^n )  \;. 
$$

The $\square $- Laplacians $\square^{(0,0)}_{\varphi }$ and $\square^{(0,1)}_{\varphi }$ are defined by
\begin{equation}\label{defsquare}
\begin{array}{ll}
\square^{(0,0)}_{\varphi } &= \ovli D^*_1 \ovli D_1 \ \ , \\
\square^{(0,1)}_{\varphi }& = \ovli D_1 \ovli D^*_1+\ovli D^*_2 \ovli
D_2\;.
\end{array}
\end{equation}

It follows that for $g=\sum_{j=1}^n g_j\, d\ovli z_j $ we have
that $\square^{(0,1)}_{\varphi }g$ equals 
$$\sum_{k=1}^n \left[ \sum_{j=1}^n \left(
2\frac{\partial^2 \varphi}{\partial z_j \partial \ovli z_k}\, g_j
-\frac{\partial^2 \varphi}{\partial z_j \partial \ovli z_j}\, g_k
-\frac{\partial^2 g_k}{\partial z_j \partial \ovli z_j}+
\frac{\partial g_k}{\partial \ovli z_j}\, \frac{\partial \varphi}{\partial z_j}-
\frac{\partial g_k}{\partial z_j}\, \frac{\partial \varphi}{\partial \ovli z_j}+
\frac{\partial \varphi}{\partial z_j}\, \frac{\partial \varphi}{\partial \ovli z_j}\, g_k \right) \right] \,
d\ovli z_k \,
$$
and that 
\begin{equation}\label{identity}
\square^{(0,1)}_{\varphi } = \square^{(0,0)}_{\varphi } \otimes I  + 2 M_{\varphi },
\end{equation}
where
\begin{equation}\label{levim}
M_{\varphi } = \left ( \frac{\partial^2 \varphi}{\partial z_j \partial
    \ovli z_k} \right )_{jk}\,.
\end{equation}

For $\varphi $ in $\mathcal C^2$,  it can be shown (by an extension of
a criterion of Simader \cite{Sima}) that $\square^{(0,1)}_{\varphi }$ can be extended to a densely defined self-adjoint operator on $L^2_{(0,1)}(\mathbb C^n),$ which is again denoted by 
$\square^{(0,1)}_{\varphi }$\,.\\

We can now state  a natural, rather standard, existence theorem for the canonical 
operator.
\begin{theorem}\label{thnew}~\\
Let us assume that $$0\not\in \sigma(\square^{(0,1)}_{\varphi })\;.
$$
Then, if $N_\varphi$ denotes its inverse, the operator 
$$
S_\varphi := (\ovli D_1)^* N_\varphi\;,
$$
is continuous from $L^2_{(0,1)}(\mathbb C^n)$ into $L^2(\mathbb C^n)$ 
 and its restriction to $\mbox{ Ker }D_2$ gives the canonical solution
 operator $S^{can}_\varphi$, hence $S^{can}$
 via \eqref{equivscan}.
\end{theorem}
\begin{proof}~\\
We have~:
$$
\begin{array}{ll}
\langle S^*_{\varphi } S_{\varphi }v\,,\,
v\rangle&
 = \langle N_{\varphi } \ovli D_1 \ovli D^*_1 N_{\varphi }v\,,\,v
 \rangle\\
&= \langle \ovli D_1 \ovli D^*_1 N_{\varphi }v\,,\,N_{\varphi }v\rangle \\
& \le \langle \ovli D_1 \ovli D^*_1 N_{\varphi }v\,,\, N_{\varphi }v\rangle +
\langle \ovli D^*_2 \ovli D_2 N_{\varphi }v\,, \,N_{\varphi }v\rangle\\
&=\langle N_{\varphi }v\,,\,v\rangle\,.
\end{array}
$$
Hence
\begin{equation}\label{compNS}
||S_{\varphi} v||^2 = \langle S^*_{\varphi } S_{\varphi
}v\,,\,v\rangle 
\le \langle N_{\varphi
}v\,,\,v\rangle\,.
\end{equation}

\end{proof}
We also indicate that
\begin{equation}\label{formuleforsquare}
4\, \square^{(0,0)}_{\varphi } = \Delta_{\varphi}^{(0)} - \Delta\varphi,
\end{equation}
where
$$\Delta_{\varphi}^{(0)}=-\sum_{j=1}^n\left ( \left (\frac{\partial}{\partial x_j}+i\frac{\partial\varphi}{\partial y_j}\right )^2 +
\left ( \frac{\partial}{\partial y_j}-i\frac{\partial\varphi}{\partial x_j}\right )^2 \right )$$
and 
$$
\Delta\varphi = \sum_{j=1}^n\left ( \frac{\partial^2
    \varphi}{\partial x_j^2}+\frac{\partial^2\varphi}{\partial
    y_j^2}\right ).$$
\vskip 1 cm

\section{About general criteria of compact resolvent}

\vskip 1 cm
The analysis of the compactness of the  canonical solution operator to $\ovprt $ involves 
 the analysis of the compact resolvent property for Schr\"odinger
 operators with compact manifold. We recall in this section  a theorem due to Helffer-Mohamed (\cite{HeMo}) on compact resolvents of Schr\"odinger operators with magnetic fields.

We will analyze the problem for the family of operators~:
\begin{equation}
P_A = \sum_{j=1}^n  (D_{x_j}- A_j(x))^2 \;.
\end{equation}
Here  $D_{x_j}=-i \frac{\partial}{\partial x_j}$ and the magnetic potential
 $A(x)=(A_1(x), A_2 (x),\cdots, A_n(x))$ is supposed to be $C^\infty$. Under these conditions, the operator is 
essentially self-adjoint on $C_0^\infty (\mathbb R^n)$. We note also that it has the form~:
$$
P_A  = \sum_{j=1}^n X_j^2\;,
$$
with
$$
X_j = (D_{x_j} - A_j(x))\,,\, j=1,\dots,n\;.\;
$$
Note that with this choice $X_j^* = X_j$.
In particular, the magnetic field
 is recovered by observing that
$$
B_{jk} = \frac 1i [X_j, X_k] = \frac{\partial A_k}{\partial x_j} - \frac{\partial A_j}{\partial x_k} \;,\;\mbox{ for } j,k = 1,\dots,n\;.
$$
We introduce for $q\geq 1$ the quantities~:
\begin{equation}
m_q(x) = \sum_{j<k}  \sum_{|\alpha|=q-1} |\partial_x^\alpha B_{jk} (x)|\;.
\end{equation}
It is easy to reinterpret this quantity in  terms of 
commutators of the $X_j$'s.\\
Let us also introduce
\begin{equation}
m^r (x) = 1 + \sum_{q=0}^r m_q(x)\;.
\end{equation}
Then the criterion is
\begin{theorem}\label{thhemo} (\cite{HeMo})~\\
Let us assume that there exists $r$ and a constant $C$ such that
\begin{equation} \label{condhemo}
m_{r+1}(x) \leq C\; m^r(x)\;,\;\forall x\in \mathbb R^n\;,
\end{equation}
and
\begin{equation}
m^r(x) \to + \infty\;,\; \mbox{ as } |x| \to + \infty\;.
\end{equation}
Then $P_A$ has a  compact resolvent.
\end{theorem}
(see also \cite{She1} and \cite{KS} for further results in this direction.)

\vskip 0.2 cm

We will mainly apply this result for the case of real dimension $2n,$ where we will write the elements of $\mathbb R^{2n}$ in the form $(x_1, y_1,\dots , x_n, y_n)$ and for the magnetic potential
\begin{equation}
A = \left ( -\frac{\partial \varphi}{\partial y_1}, \frac{\partial \varphi}{\partial x_1}, \dots ,
-\frac{\partial \varphi}{\partial y_n}, \frac{\partial \varphi}{\partial x_n}\right) \;.
\end{equation}

\vskip 1 cm

\section{The analysis of the Laplacian and application}

\begin{theorem}\label{Th32}~\\
Let $\varphi $ be a plurisubharmonic $\mathcal C^2$ - function on $\mathbb C^n$ such that for the lowest eigenvalue $\lambda_{\varphi}$ of the Levi matrix $M_{\varphi}$ the condition
\begin{equation}\label{asslev}
\liminf_{|z|\to \infty} \lambda_{\varphi} (z) >0 \,,
\end{equation}
 is satisfied. Then the operator $\square^{(0,1)}_{\varphi }$ has a bounded inverse $N_\varphi $ on $L^2_{(0,1)}(\mathbb C^n)$.
\end{theorem}
\begin{proof}~\\
For $v = \sum_{k=1}^n v_k\, d\ovli z_k \in
\textrm{Dom\,}\square_{\varphi }^{(0,1)} $, we have by \eqref{identity},
\begin{equation}\label{comparisona}
\langle \square_{\varphi }^{(0,1)} v\,,\, v\rangle 
\ge 2\langle M_{\varphi }v\,,\,v\rangle\; .
\end{equation}
Using Persson's Theorem (see for instance \cite{Ag}),  we now conclude
from Assumption \eqref{asslev} that the bottom of the essential
spectrum of $\square_{\varphi }^{(0,1)}$ is strictly positive.
Using the spectral theorem for selfadjoint operators,  we conclude that 
$\square_{\varphi }^{(0,1)} $ is bijective if $\square_{\varphi }^{(0,1)}$ is injective
(see for instance \cite{W} , (8.17)). 
In order to show that $\square_{\varphi }^{(0,1)}$ is injective we consider the inequality
\begin{equation}\label{minor1}
\langle \square_{\varphi }^{(0,1)} v\,,\, v\rangle \ge \int_{\mathbb
  C^n}\sum_{k=1}^n \lambda_{\varphi} (z) \vert v_k(z)\vert ^2 \,
d\lambda (z).
\end{equation}
We recall that $\lambda_\varphi \geq 0$. 
If $ \square_{\varphi }^{(0,1)} v=0,$ \eqref{minor1} together with 
 Assumption \eqref{asslev} implies that $\lambda_\varphi$ is non zero
 at $\infty$, hence 
$v=0$ on a non-empty open set. Therefore by the uniqueness result of
 Kazdan (\cite{Kaz}) it follows that $v=0$ everywhere and that $
 \square_{\varphi }^{(0,1)} $ is injective and therefore also
 surjective and has a bounded inverse $N_{\varphi }$\,. Hence we can
 apply Theorem \ref{thnew}.
\end{proof}

\begin{theorem}\label{Th33}~\\
Let $\varphi $ be a plurisubharmonic $\mathcal C^2$ - function on
$\mathbb C^n$ such that
 \begin{equation}\label{assl2}
\lim_{|z|\rightarrow \infty}\lambda_\varphi(z)  = +\infty\;.
\end{equation}
 Then the canonical solution  operator to $\ovli \partial$  $S^{can}$ is compact.
\end{theorem}
 \begin{proof}~\\
By Theorem \ref{thnew} and \eqref{equivscan}, it is sufficient to show that $S_\varphi$ is compact from
$L^2_{(0,1)}(\mathbb C^n)$ into $L^2(\mathbb C^n)$. 
Using \eqref{minor1} and \eqref{assl2}, it follows that  $\square_{\varphi }^{(0,1)}$ has compact resolvent 
(see for instance \cite{AHS} or \cite{I}) and we have also shown 
in Theorem \eqref{Th32},  that 
$\square_{\varphi }^{(0,1)} $ was bijective. The operator  $N_{\varphi
}$ 
is consequently a compact self-adjoint operator on 
$L^2_{(0,1)}(\mathbb C^n).$

The operator $S_{\varphi}=\ovli D^*_1 N_{\varphi }$ is the canonical solution operator to $\ovli D_1v=g\,$.
Now if $N_\varphi$ is compact, it is standard that $ N_\varphi^\frac
12$ is compact. It is then easy to show from \eqref{compNS}
that
$S_{\varphi }$ is compact.

 \end{proof}

\begin{rem}~\\
Theorem \ref{Th33} can be applied for instance in the case when the weight function is of the form 
$$
\varphi (z)=\left (\sum_{j=1}^n |z_j|^2 \right )^m\,, 
$$
for some integer $m>1.$ This is strongly related to examples given by M. Derridj for the analysis of the regularity of 
$\square_b ,$ as discussed in the book \cite{HeNo} (Chap. V.2).
\end{rem}

\begin{rem}~\\
Theorem \ref{Th33} should be compared with the corresponding estimate in \cite{H} (4.4.1.), which is of the form
$$\int_{\mathbb C^n} | u(z) |^2 \, e^{-2\varphi (z)}\,d\lambda (z) \le
\, \int_{\mathbb C^n} | \ovli \partial u(z) |^2 \, \,
\frac{e^{-2\varphi (z)}}{\lambda_{\varphi} (z)}\, \,d\lambda (z),$$
for all $u$ in the domain of $\ovli \partial$ orthogonal to 
$\ker \ovli \partial$.

In addition we note that the last inequality is similar to
 a Brascamp-Lieb inequality as analyzed
 by Witten-Laplacians techniques (see for example \cite{He3} and the
 references therein including the generalization obtained by \cite{Jo}). 

 If  $0$ is not in the spectrum of $\square_{\varphi}^{(0,1)}$, then
 we have 
\begin{equation}\label{BraLie}
\int_{\mathbb C^n} | u(z) |^2 \, e^{-2\varphi (z)}\,d\lambda (z) \le
\frac 12  \langle M_{\varphi}^{-1}  \ovli \partial u\;,\;  \ovli \partial u
\rangle_{L^2_{(0,1)}(\mathbb C^n,\varphi)}
,
\end{equation}
for all $u$ in the domain of $\ovli \partial$ orthogonal to 
$\ker \ovli \partial$.\\
Let us give the very short proof.
By Ruelle's Lemma \cite{Ru},  we immediately deduce 
 from \eqref{comparisona} 
that
$$ N_{\varphi} \leq \frac 12 M_{\varphi}^{-1}\;.$$
Now, with $v=u \exp (- \varphi)$ and $g =\overline{D_1} v =\exp (-\varphi)\, \overline{\partial} u$, we obtain~:
$$
|| v||^2 = \langle v\,,\, S_\varphi g\rangle = \langle g \,,\,
N_{\varphi}
g\rangle \leq \frac 12 \langle  M_{\varphi}^{-1} g\,,\, g \rangle\;,
$$
 where all the norms and scalar products are in $L^2$ with the
 Lebesgue measure. This gives \eqref{BraLie}.

This implies in particular H\"ormander's statement above, but not
 Shigekawa's result below. 
\end{rem}

\begin{rem}~\\
In this connection it is also interesting to mention a result of Shigekawa (\cite{Shi}) stating that the space $A^2(\mathbb C^n , \varphi )$ is of infinite dimension if the lowest eigenvalue $\lambda_{\varphi} (z)$ of $M_{\varphi }$ satisfies the condition 
$$
\lim_{\vert z \vert \to \infty}\vert z \vert^2\,\lambda_{\varphi}
(z)  = \infty\, .
$$
This condition implies that $0$ is an eigenvalue of infinite multiplicity for a  Pauli operator of the form
$$\widetilde{H} (a) = \sum_{j=1}^{2n} (-i\partial_j - a_j(x))^2 + \sum_{j,k=1}^{2n} \frac{i}{2}\, b_{jk}(x)\gamma^j\gamma^k,$$
acting on $L^2(\mathbb R^{2n})\otimes \mathbb C^r,$ where $b_{jk}=\partial_j a_k-\partial_k a_j,$ where $r=2^n$ 
and where the $\gamma^j$'s are the $r\times r$ Dirac matrices satisfying $\gamma^j \gamma^k + \gamma^k \gamma^j = 2\delta^{jk},$ ($\delta^{jk}$ being the Kronecker delta)
(see \cite{Shi}).

Shigekawa also analyzes the link between $\widetilde{H} (a)$ and the complex Witten Laplacian by comparing the essential spectra of these operators.
\end{rem}

\vskip 0.2 cm
Finally we prove a variant of Theorem \ref{Th33} using the results
from \cite{HeMo}, together with ideas of M.~Derridj (see \cite{HeNo}
and references therein).
\vskip 0.2 cm

\begin{theorem}\label{Th37}~\\
If $\varphi $ is a plurisubharmonic $\mathcal C^2$ - function on $\mathbb C^n$ and suppose that there exists a number $t\in (0,1/4)$ and a compact set $K$ in $\mathbb C^n$ such that for the Levi matrix $M_{\varphi } $ the estimate
$$M_{\varphi } \ge t \Delta \varphi \otimes I$$
holds outside of $K$ and that $\lambda_{\varphi}$ does not vanish identically. Assume that $\Delta_\varphi $ has compact resolvent.
Then the canonical solution operator $S$ operator to $\ovprt $
is compact.
\end{theorem}
\begin{proof}~\\
Using \eqref{identity}, we have~:
 \begin{equation}
\square_\varphi ^{(0,1)}\geq  (\square_\varphi ^{(0,0)} + 2t \Delta \varphi)\otimes I\;,
\end{equation} 
outside the compact set $K$.\\

By formula \eqref{formuleforsquare},  we are then reduced to the analysis of the compactness of the resolvent 
 of
$$ 
\frac{1}{4}\, \Delta_\varphi^{(0)}  + (2t-1/4)  \Delta \varphi\;.
$$
 which is reduced, observing that for some constant $C_t >0$
 we have
$$ \frac{1}{C_t} \Delta_\varphi^{(0)}\leq  \Delta_\varphi^{(0)}  + (8t-1)
\Delta \varphi \leq C_t  \Delta_\varphi^{(0)}\;,
$$
to the same question for
$\Delta_\varphi^{(0)}$.\\

\end{proof}

We now complete the discussion by saying under which condition
$\Delta_\varphi^{(0)} $ has compact resolvent.
 This has been done already  in detail when $n=1$. One can of course
 use the criterion of Helffer-Mohamed recalled  in the previous section
 (or some of the improvements obtained later).\\
Actually, a complementary  result can be obtained by generalizing our
analysis in $\mathbb C$. We observe indeed in the same way as in the
case of $\mathbb C$, 
 that $\Delta_\varphi^{(0)} $ has compact resolvent if $-\Delta + \Delta
 \varphi$
 has compact resolvent. \\ This is then the case
 if we have the conditions that $\Delta \varphi\in A_\infty$
 and if
$$
\liminf _{|z|\rightarrow \infty} \int_{\Pi_{j=1}^n B(z_j,1)} (\Delta
\varphi(y))\,d\lambda  = +\infty\;.
$$

\vskip 1 cm
\section{The case of decoupled weights}~\\

\vskip 0.1 cm

Here we consider weights $\varphi $ of the form
$$\varphi (z_1, \dots , z_n)=\sum_{j=1}^n \varphi_j(z_j),$$
where the functions $\varphi_j $ are  $C^\infty$ on $\mathbb C$.
~\\

\vskip 0.6 cm

\subsection{About Dirac and Pauli operators}~\\

\vskip 0.3 cm

In this case an interesting connection to Dirac and Pauli operators is of importance (see \cite{CFKS}, \cite{Er}, \cite{HNW}, \cite{Roz}, \cite{Tha}).
Let us first consider the real two dimensional case. The Dirac operator $\mathbb D$ is defined by
$$
\mathbb D =  \sigma_1 \left ( \frac{1}{i} \, \partial_{x_1}- A_1(x,y)
\right ) +  \sigma_2 \left ( \frac{1}{i} \, \partial_{x_2}- A_2(x,y) \right ),$$
where
$$
\sigma_1=\left(
\begin{array}{cc}
 0&1\\1&0
\end{array}\right)\;,\;\sigma_2=\left(
\begin{array}{cc}
 0&-i\\i&0
\end{array}\right)\;.
$$
It turns out that the square of $\mathbb D$ is diagonal with the Pauli operators $P_{\pm }$ on the diagonal:
$$\mathbb D^2 = \left(
\begin{array}{cc}
 P_-&0\\0&P_+
\end{array}\right),$$
where 
$$P_{\pm}= 
 \left ( \frac{1}{i} \, \partial_{x_1}- A_1(x,y) \right )^2
 +  \left ( \frac{1}{i} \, \partial_{x_2}- A_2(x,y) \right )^2 \pm
 B(x,y)\,.
$$
Using the computation done in \eqref{formuleS}, we get, having in mind that
 $\mathcal S = \square^{(0,0)}_{\varphi }$, 
$$
4\, \square^{(0,0)}_{\varphi } = P_- \;.
$$

It is proved in \cite{HNW} (Theorem 1.3) that at least one of the
operators $P_{\pm}$ has non compact resolvent if $\varphi$
satisfies in $\mathbb C$ 
 the following condition ($H_r$)~:\\
There exists a sequence of disjoint balls $B_n$ of radius $\geq 1$
 such that \eqref{condhemo} is satisfied in the union of these balls.\\

This is in particular the case when the magnetic potentials are
polynomials.\\

Note also the interesting independent result (cf \cite{CFKS}) that the
spectra of $P_+$ and $P_-$ coincide
 except at $0$. So if $P_+$ has compact resolvent then $P_-$ has its
 essential spectrum
 reduced to $\{0\}$.\\

 \vskip 0.3 cm
 
\subsection{Main results and proofs}~\\

\vskip 0.3 cm

Our main theorem in this section is the following
\begin{theorem}\label{apllHNW}~\\
Let $n\geq 2$ and let  $\varphi$ be  a decoupled weight such that 
 there exists $j$ such that  $\varphi_j$
 satisfies for some $r_j>0$ the condition $(H_{r_j})$, then
 $\square_{\varphi}^{(0,1)}$ has a non compact resolvent.
\end{theorem}
\begin{proof}~\\
As observed in  \cite{Has3}, a simple computation shows that
 for the decoupled weights 
$$\varphi (z_1, \dots , z_n)=\sum_{j=1}^n \varphi_j(z_j)$$
the operator $\square^{(0,1)}_{\varphi } $ becomes diagonal,
 each component on the diagonal being
\begin{equation}
\mathcal S_k= \square^{(0,0)}_{\varphi }
 + 2 \ \frac{\partial^ 2 \varphi_k}{\partial z_k \partial \ovli z_k}.
\end{equation}
Then the result is based on the following proposition.
\begin{prop}~\\
Let $n\geq 2.$ Under the  assumptions of the theorem  on the weight function $\varphi$, there always exists a $k$ such that $\mathcal S_k$ is not with compact resolvent.
\end{prop}
We observe that $\mathcal S_k$  can be rewritten in the form 
$$
4\mathcal S_k = \sum_{j\neq k} P_-^{(j)} + P_+^{(k)}\,,
$$
where each operator $P_{\pm}^{(\ell)}$ is the previously 
analyzed Pauli operator in variables the $(x_{\ell}, y_{\ell}).$
The result is then obtained from the results by Helffer-Nourrigat-Wang
recalled in the previous subsection.
\end{proof}

\vskip 0.3 cm

\begin{rem}~\\
It is also easy to see that the kernel of $P_-^{(\ell )}$ contains all $L^2$-distributions of the form 
$$
f(z_{\ell})
\exp (-\varphi_{\ell}(z_{\ell})),
$$
where $f$ is holomorphic and $z_{\ell }=x_{\ell }+iy_{\ell }\,$.

Hence $\mathcal S_k$ has non-compact resolvent,
 as soon as the space $A^2(\mathbb C, \varphi_\ell)$ is of infinite
 dimension
 for some $\ell\neq k$. This can be combined with Shigekawa's result,
 see
 also the next propositions.
\end{rem}
~\\
\vskip 0.3 cm

\subsection{On a result of G.~Schneider}~\\

\vskip 0.3 cm

In the case of decoupled weights, one can extend a remark of
G. Schneider (\cite{Sch}) who was considering  the case when
 $\varphi_j(z_j)= |z_j|^{2m}$
 for $m>1$,  to show 
that the canonical solution operator to $\ovprt $ fails to be compact
even on the space $A^2_{(0,1)}(\mathbb C^n, \varphi )$ of
$(0,1)$-forms with holomorphic coefficients. \\

\begin{prop}\label{prop6.4}~\\
 Suppose that $n\geq 2$ and that  
there exists $\ell$ such that
 $A^2(\mathbb C,\varphi_\ell)$ is infinite dimensional. Suppose also that
 $ 1 \in L^2(\mathbb C,\varphi_j)$ for all $j$
 and that there exists $k\neq \ell$ such that
 $\frac{ \partial^2\varphi_k}{\partial z_k \partial
 \ovli z_k}   \in L^2(\mathbb C,\varphi_k)$. Then
 $\mathcal S_k$ has non compact resolvent. In particular,
  $\square_\varphi^{(0,1)}$ has non compact resolvent.
\end{prop}

\begin{proof}~\\
Let $f_\nu$ an infinite orthonormal system in $A^2(\mathbb
C,\varphi_\ell)$. 
For the functions
$$
u_\nu(z)= f_\nu (z_\ell)\,\, \exp (-\varphi (z))$$
we have by \eqref{defsquare} 
$$
\square^{(0,0)}_{\varphi }u_\nu = \ovli D^*_1 \ovli D_1 u_\nu =0\;,
$$
for all $\nu =1,2,\dots $ and by \eqref{identity} 
$$
\square^{(0,1)}_{\varphi } \left ( u_\nu\, d\ovli z_{k }\right ) 
=\left(\mathcal S_\ell u_\nu\right) d\ovli z_k= \left( 2 \ \frac{\partial^ 2 \varphi_{k}}{\partial z_{k } \partial \ovli z_{k}}
u_\nu  \right) \, d\ovli z_{k}\,.
$$

Hence, the sequence 
$$
\langle \square^{(0,1)}_{\varphi } \left ( u_\nu\, d\ovli z_{k
  }\right )\, , \,\left ( u_\nu\, d\ovli z_{k }\right )\rangle =\langle
  \mathcal S_k\,u_\nu\;,\; u_\nu
\rangle $$
is bounded 
and, by the assumption that the functions  $z_\ell \mapsto
  f_\nu(z_\ell) \, \exp (-\varphi_\ell(z_\ell))$ form an orthonormal system, we
  get the statement.
\end{proof}

Using a similar argument we get the following extension of a result of G. Schneider \cite{Sch} (see also \cite{Kr}).

\begin{prop}~\\
Suppose that $n\geq 2$ and that  
there exists $\ell$ such that
 $A^2(\mathbb C,\varphi_\ell)$ is infinite dimensional. Suppose also  that
 $1\in L^2(\mathbb C,\varphi_j)$ for all $j$. Suppose finally 
 that for
some
 $k\neq \ell$, $\ovli z_k \in L^2(\mathbb C,\varphi_k)$. Then the
canonical solution operator to $\ovli \partial$ fails to be compact
even on the space $A^2_{(0,1)}(\mathbb C^n,\varphi)$.
\end{prop}

\begin{proof}~\\
Let $P_k$ denote the Bergman projection from $L^2(\mathbb
C,\varphi_k)$
 onto $A^2(\mathbb C, \varphi_k)$.
It is clear that the function $(\ovli z_k -P_k \ovli z_k)$ is not zero. 
With the notations of the preceding proof, the family 
$$h_\nu :=f_\nu
(z_\ell )
(\ovli z_k -P_k \ovli z_k)
$$ 
is an orthogonal family in  $A^2(\mathbb
C^n,\varphi)^{\perp}$, which 
 satisfies $\ovli \partial h_\nu = f_\nu (z_\ell) d\ovli z_k$.

Hence $ \left ( \ovprt h_{\nu} \right )_{\nu }$ constitutes a bounded sequence in $A_{(0,1)}^2 (\mathbb C^n, \varphi ),$
and this  implies the result.
\end{proof}

\end{document}